\documentclass[conference]{IEEEtran}
\usepackage{cite}
\usepackage{cite}
\usepackage[pdftex]{graphicx}
\usepackage{amsmath}
\usepackage{fixltx2e}
\usepackage{pgfplots}
\usetikzlibrary{matrix,positioning}
\usepackage{tikz}
\usepackage{booktabs, threeparttable}
\usepackage{multirow}
\usepackage{array}

\ifCLASSINFOpdf
\else
\fi
\begin{document}
%
\title{A Centralized Power Control and Management Method for Grid-Connected Photovoltaic (PV)-Battery Systems}

\author{

\IEEEauthorblockN{Zhehan Yi}
\IEEEauthorblockA{Department of Electrical and\\Computer Engineering\\
	The George Washington University\\
	Washington, DC 20052\\
	Email: zhehanyi@gwu.edu}
\and
\IEEEauthorblockN{Wanxin Dong}
\IEEEauthorblockA{Department of Electrical and\\Computer Engineering\\
	The George Washington University\\
	Washington, DC 20052\\
	Email: wanxindong@gwu.edu}
\and
\IEEEauthorblockN{Amir H. Etemadi}
\IEEEauthorblockA{Department of Electrical and\\Computer Engineering\\
	The George Washington University\\
	Washington, DC 20052\\
	Email: etemadi@gwu.edu}
}


%


\maketitle

\begin{abstract}
Distributed Generation (DG) is an effective way of integrating renewable energy sources to conventional power grid, which improves the reliability and efficiency of power systems. Photovoltaic (PV) systems are ideal DGs thanks to their attractive benefits, such as availability of solar energy and low installation costs. Battery groups are used in PV systems to balance the power flows and eliminate power fluctuations due to change of operating condition, e.g., irradiance and temperature variation. In an attempt to effectively manage the power flows, this paper presents a novel power control and management system for grid-connected PV-Battery systems. The proposed system realizes the maximum power point tracking (MPPT) of the PV panels, stabilization of the DC bus voltage for load plug-and-play access, balance among the power flows, and quick response of both active and reactive power demands.      
\end{abstract}


%
\IEEEpeerreviewmaketitle

\section{Introduction}
Contributions of renewable energy to power generation have been increasing exponentially in the past decades, due largely to the fossil energy crisis, increasing electricity demands, and environment degradation. Power grids in the near future is expected to be equipped with great penetration of renewable energy with effective supply-demand management systems for highly reliable and economical operations \cite{yi_tsg2}. Distributed generation (DG), through which electric power is generated on-site instead of centrally, is providing a powerful solution of integrating renewable power generators to the conventional utility grids. Among numerous renewable generations, solar photovoltaic (PV) system is one of the most attractive renewable power system because of its various benefits, such as flexibility of scales and low installation costs. Moreover, as the price continues to decline, global PV installation capacity is expected to increase in recent years \cite{yi_tsg}. As the penetration of solar power expands, future grid-connected PV systems are required to provide more reliable power. Otherwise, power fluctuations in PV systems will bring certain reliability issues to the utility grids and to electricity users.    

PV output power oscillates frequently during a day as the operating environment (temperature or solar irradiance) changes \cite{yi_tsg, yi_tie, yi_pes2016}. Therefore, to maintain a stable output, battery storages are necessary on the DC side of a PV system to compensate the differences caused by change of operating condition, for example, cloud shading over the PV panels. Excess power can also be stored in the batteries for later demands, or to trade back to the utility in the future. Additionally, there are usually loads on the DC bus, making the system a DC microgrid. Therefore, the control schemes for PV-battery systems should be able to stabilize the power supply to the DC loads, balance the power flows inside the DC side, and effectively manage the power communications with the grid. 

A number of control schemes for PV-battery units have been proposed in the literature. An autonomous control strategy is
presented in \cite{literature5} for PV-battery systems by droop control. The battery group is charged through the AC bus, which increases the costs for the charger inverter. Furthermore, this strategy only works for islanded but not for grid-connected PV systems. Other islanded PV control systems are also proposed in references \cite{literature7} and \cite{literature11}. Again, these methods are not applicable to grid-connected PV systems, which are widely used in the industry. A hierarchical control scheme is designed in \cite{literature6} for a grid-connected multi-source PV systems, which is primarily used for self-feeding buildings equipped with PV arrays and battery storages. This system requires complicated supervision algorithm using Petri nets (PNs), and the DC bus voltage oscillations is not well eliminated. Literature \cite{literature8} presents an optimal charging/discharging method in PV-battery system that reduces the line loss of distribution systems. However, this method only schedules the battery changing and discharging but does not comprehensively manage the power flows of the entire system. There are also other control strategies introduced in \cite{literature9,literature10}, which mainly focus on the optimization of operational costs. Nevertheless, these paper does not introduce the detailed control methods.  
 
In an attempt to address the aforementioned issues, this paper proposes a power management and control system for grid-connected PV-battery power systems, which balances the power flows flexibly and maintains a reliable power supply to demands in different circumstances. Control methods in the system are designed to achieve a flexible but reliable power output to fulfill demands from the utility grid, and the loads on the DC bus, by intelligently managing the charging/discharging processes of the battery and switching the power generating control modes for the PV arrays. In the following contents: Section II briefly introduces a typical grid-connected PV-battery system, followed by the proposed power control and management system with detailed control schemes for each part; case studies are carried out in Section III to verify its performance; Section IV presents the conclusion of this paper.       
     

\section{Grid-Connected PV-Battery Systems and the Proposed Power Control and Management system}

\subsection{Typical Configuration of Grid-Connected PV-Battery System}

\begin{figure}
	\centering
	\includegraphics [scale=0.9] {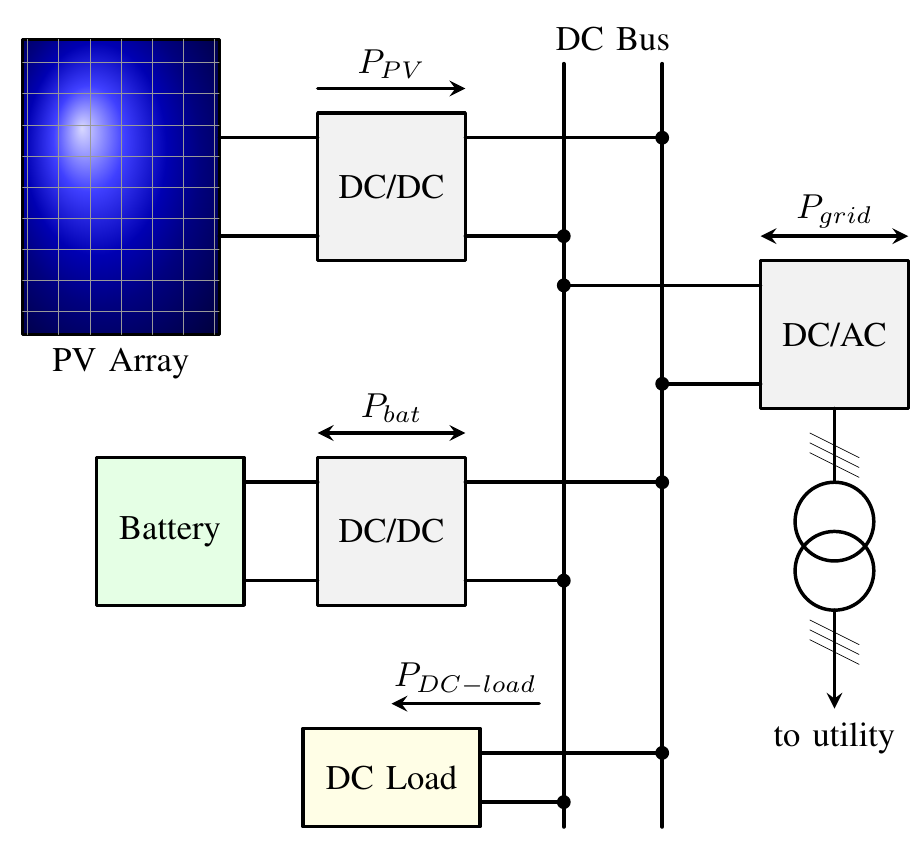}
	\caption{A typical grid-connected PV-Battery system.}\label{system}
\end{figure}

A typical PV-battery system is illustrated in Fig. \ref{system}, which consists of a PV array, battery storages, DC/DC converters that connect the PV array and the battery to the DC bus, DC load, a DC/AC inverter and a transformer to bridge the DC microgrid to the utility grid \cite{yi_tsg2, yi_phd}. The power generated by the PV array is a function of irradiance during a day, and due to the non-linear characteristics of PV, for different irradiance level, there is a operating voltage V\textsubscript{MPP} where the PV array is extracted the maximum power. Therefore, maximum power point tracking (MPPT) algorithms are implemented in PV systems, usually through a DC/DC converter, to optimize the power generation \cite{ref12,ref13,ref14}. The charging or discharging process of the battery storage is controlled by another DC/DC converter. DC loads, on the other hand, can be plugged in the DC microgrid directly or through converters. Therefore, if the control strategies are able to stabilize the DC bus voltage while balancing the power flow, converters can be omitted for load access, making the system more convenient and economical.

\begin{figure*}
	\centering
	\includegraphics [scale=0.55] {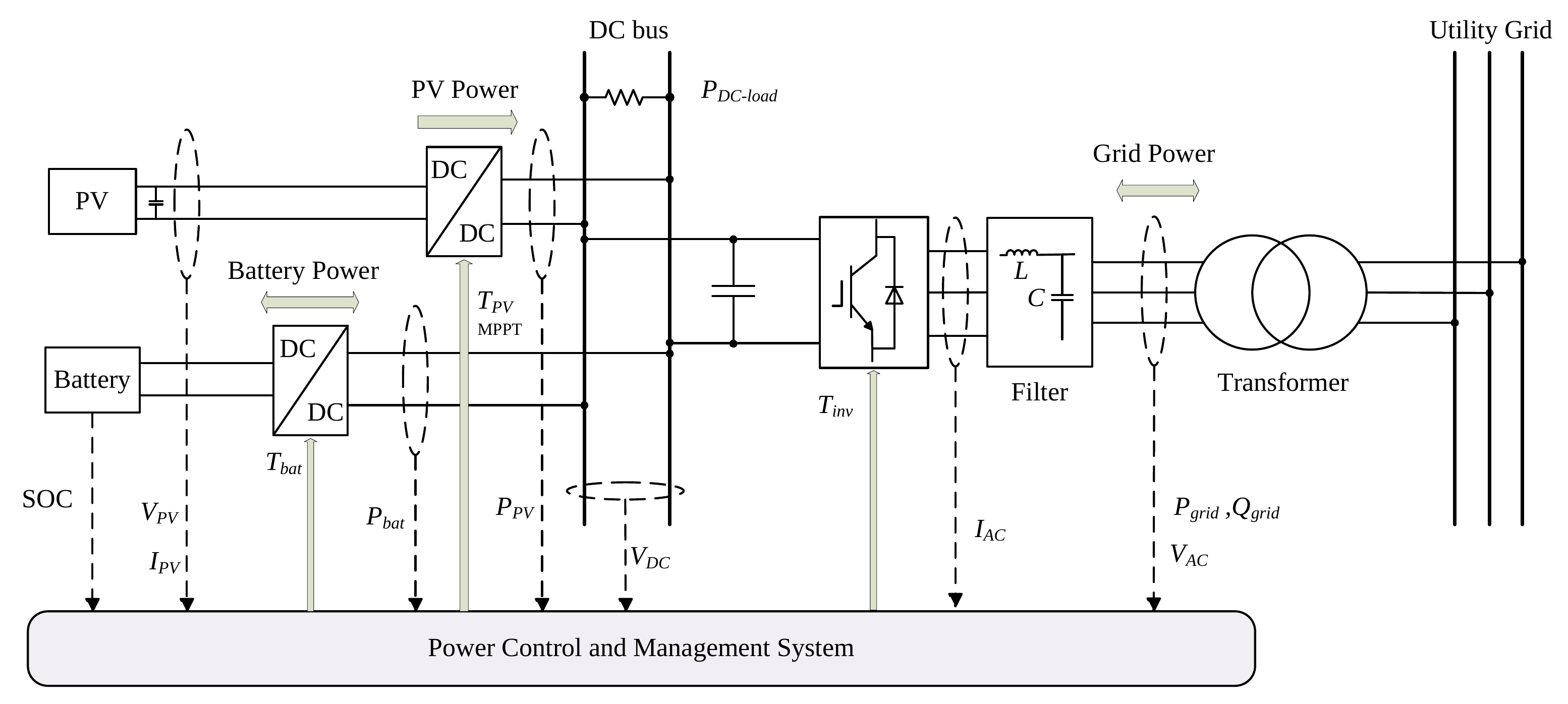}
	\caption{The proposed power control and management system for grid-connected PV-battery power systems.}\label{system2}
\end{figure*}

\subsection{The Power Control and Management System}

While some power flows in the PV-battery systems have to be unidirectional, e.g., power flows from the PV to the DC bus ($P_{PV}$) and power consumed by the DC load ($P_{DC-load}$), some should be bidirectional, e.g., power charging or discharging the battery ($P_{bat}$) and power communicates with the utility grid ($P_{grid}$). Therefore, to keep the power balanced in the system, the following equation should always be fulfilled:

\begin{equation}
P_{PV}+P_{bat}=P_{grid}+P_{DC-load}+P_{loss}
\label{eq1}
\end{equation}

\noindent where $P_{loss}$ is the power loss in the power converters, transformer, and transmission lines, which is usually negligible. $P_{bat}>0$ and $P_{bat}<0$ indicates discharging and charging mode of the battery, respectively; $P_{grid}>0$ means the power is being transferred from the DC microgrid to the utility grid, and vice versa. A power control and management system is designed (Fig. \ref{system2}), which supervises the status of each generation and load in the system, and, depending on the situations, determines the references for the PV power, DC bus voltage, battery charging/discharging power, and the active and reactive power though the inverter. Detailed control schemes are elaborated as follows.

\subsubsection{PV Generation Control}
\begin{figure}
	\centering
	\includegraphics [scale=0.55] {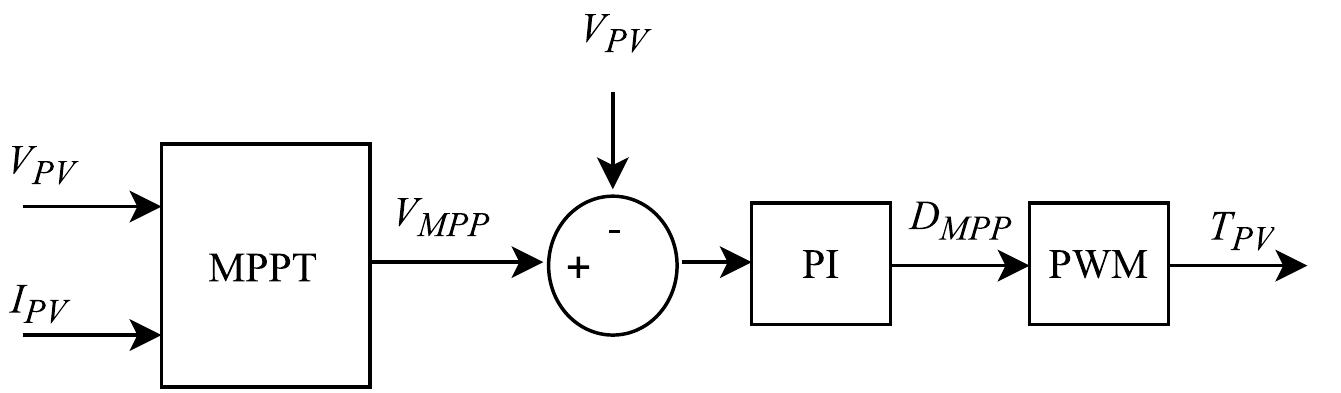}
	\caption{MPPT control of the PV array.}\label{MPPT}
\end{figure}
The power control for the PV array can be switched between MPPT control mode and power reference control mode, depending on the SOC of battery, DC load demand $P_{DC-load}$ and the grid-requested power $P_{grid}$. When the battery is not fully charged ($\textup{SOC} < 95\%$), $P_{DC-load}$ is met, and the utility grid is requesting the DC microgrid system to provide maximum power as it can, the PV array will be controlled under MPPT mode. On the other hand, if the battery is fully changed  ($\textup{SOC} \geq 95\%$), DC load demand ($P_{DC-load}$) is fulfilled, and the grid is not able to consumed the excess PV power, MPPT will turn off, and the PV will be switched to power reference control mode, where the reference is given by the following equation.
\begin{equation}
P_{PV-ref}= P_{grid}+P_{DC-load} 
\end{equation}
\noindent When working under MPPT mode, the voltage and current of the PV array, $V_{PV}$ and $I_{PV}$, are extracted and fed into the management system to obtain a voltage reference $V_{MPP}$, and generate a gating signal $T_{PV}$ for the switching control of the DC/DC (boost) converter. Incremental conductance (IncCond) MPPT, one of the most well-known MPPT algorithm \cite{ref16}, is employed in this research. The control scheme is presented in Fig. \ref{MPPT}.

\subsubsection{Battery Charging/Discharging Control}
 
\begin{figure}
	\centering
	\includegraphics [scale=0.55] {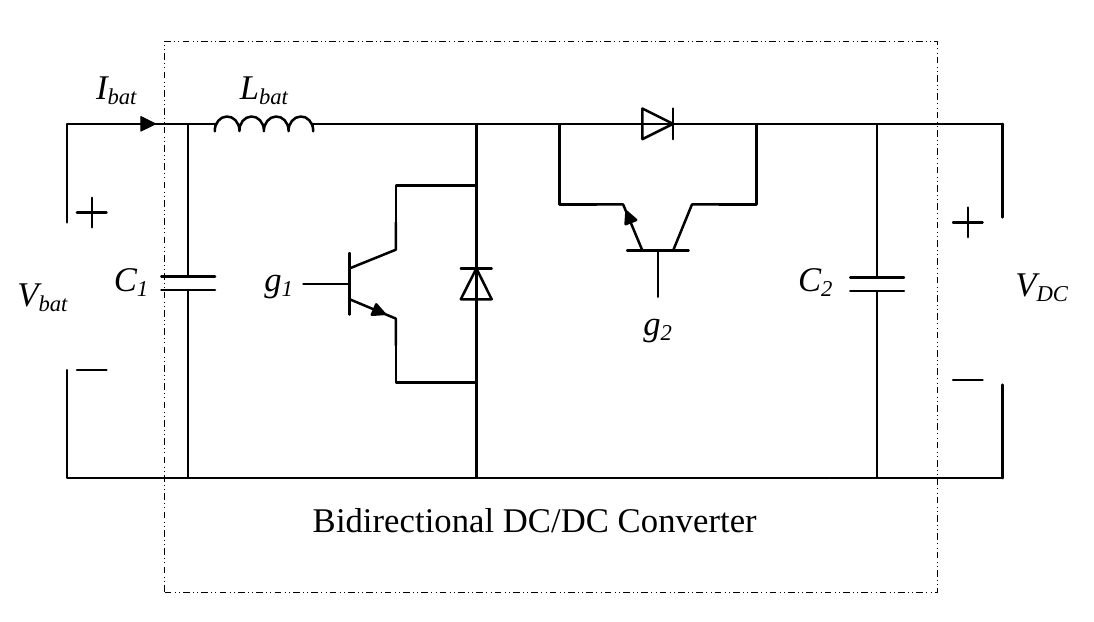}\\
	(a)\\
	\vspace{3mm}
	\includegraphics [scale=0.55] {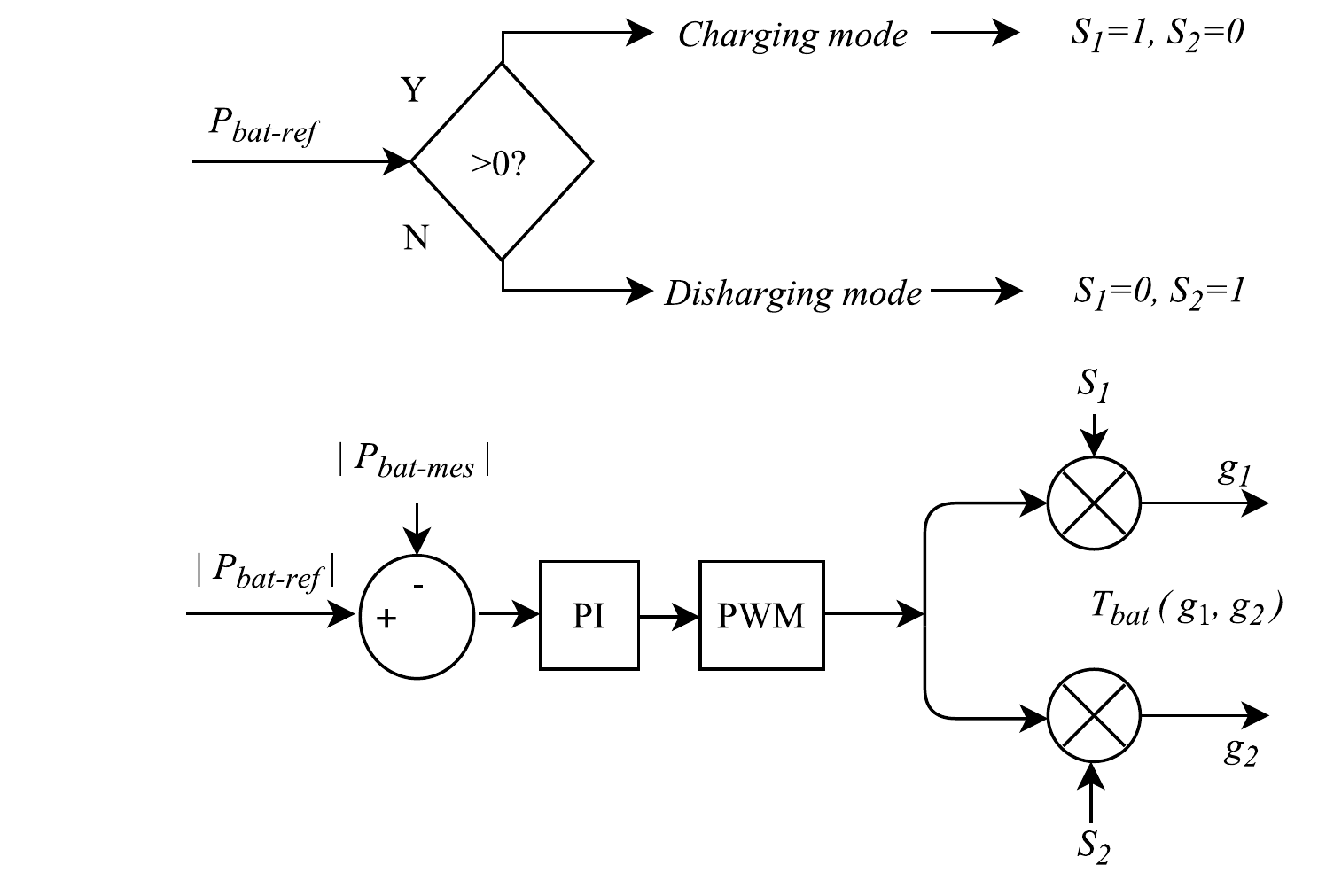}\\
	(b)
	\caption{(a) Bidirectional DC/DC converter for battery charging/discharging; (b) Battery charging/discharging control scheme.}\label{bat}
\end{figure}

The operating mode of battery, i.e., charging/discharging or the sign of $P_{bat}$, is not only subject to equation (\ref{eq1}), but also depending on the state of change (SOC) of the battery. Namely, there are maximum and minimum limits for SOC, which is set to be 95\% and 20\% respectively in this scheme, to eliminate the degradations and extend the life cycle of the battery. A bidirectional DC/DC converter is used to control the charging and discharging of battery (Fig. \ref{bat} (a)), and the control scheme is illustrated in Fig. \ref{bat} (b), where $P_{bat-mes}$ and $P_{bat-ref}$ is the measured and desired power flowing in the converter, respectively.

\subsubsection{Inverter Control} 
\begin{figure}
	\centering
	\includegraphics [scale=0.55] {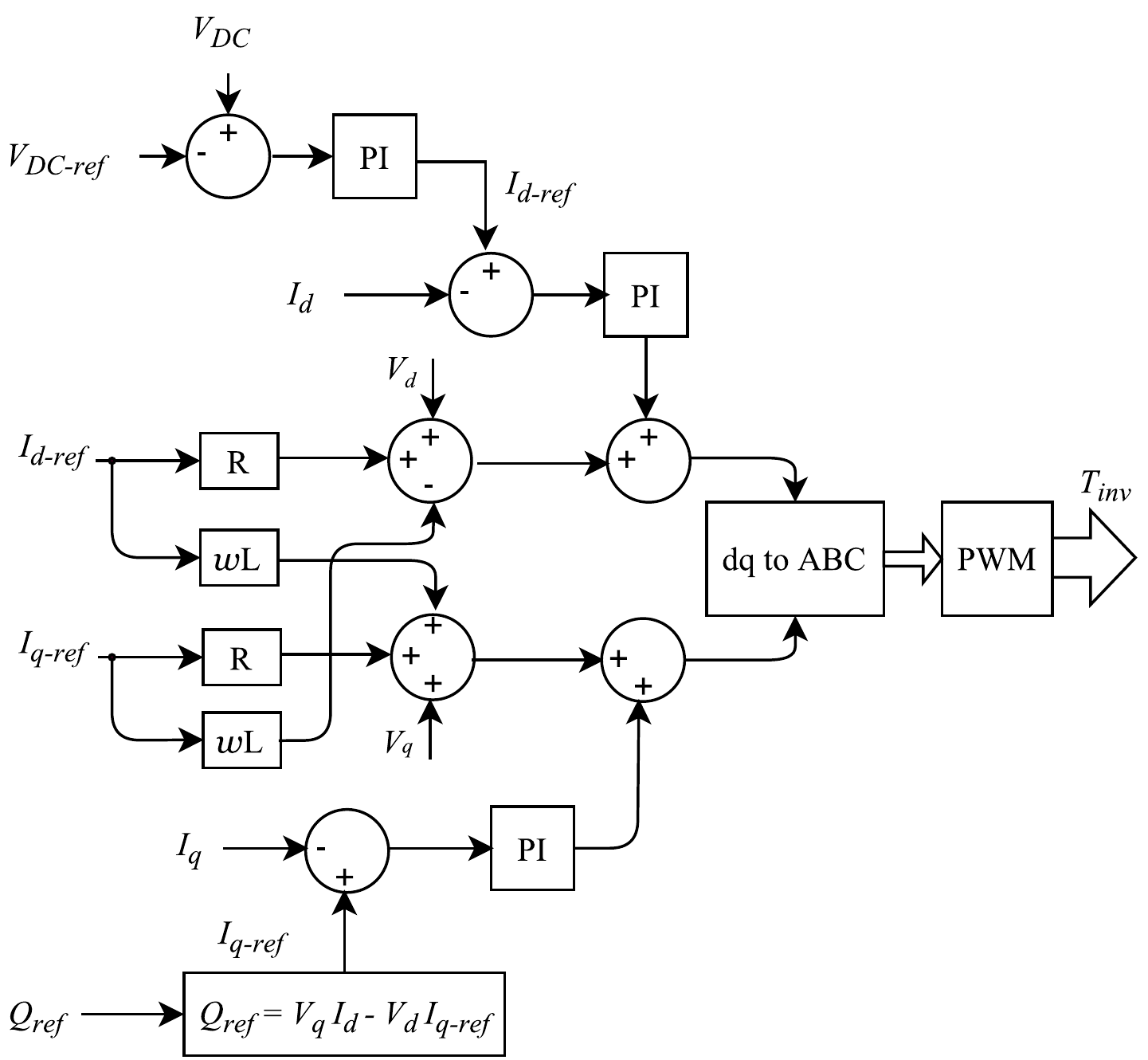}
	\caption{Inverter control scheme.}\label{inverter}
\end{figure}
An inverter is necessary to convert the DC power to AC power, and connect the DC microgrid to the utility grid via the point of common coupling (PCC). Fig. \ref{inverter} illustrates the control scheme for the inverter, which aims at stabilizing the DC bus voltage, $V_{DC}$, and controlling the reactive power, $Q$, flowing through the inverter. $I_{d}$ and $I_{q}$ are the obtained from the three-phase AC current $I_{AC}$ (Fig. \ref{system2}) by Clark transform. The active power flows in the inverter is controlled according equation (\ref{eq1}), i.e., controlling $P_{PV}$ and $P_{bat}$ so that $P_{grid}$ is controlled.     

\section{Case Studies}

In order to verify the proposed control strategies discussed in Section II, case studies are carried in this section using the PSCAD software package. A grid-connected PV-battery is set up using the configuration in Fig. \ref{system2}, where the parameters are listed in Table. \ref{parameters}. According to various statuses of the battery SOC, PV generation ($P_{PV}$), DC load demand, and utility demand ($P_{grid}$), multiple cases are simulated and the results are presented as follows.

\begin{table}
	\centering
	\caption{PV-Battery System Parameters}\label{parameters}
	
	\begin{tabular}{cc}
		\toprule[1.5pt]
		\bfseries{Parameters}  & \bfseries{Values} \\
		\midrule
		PV Maximum Power &  165\,kW \\
		DC Load in Case 1 $\sim$ Case 4 & 50\,kW\\
		DC Load in Case 5 & 190\,kW\\		
		DC Bus Voltage ($V_{DC}$) & 450\,V\\
		AC Bus Voltage (line to line) & 208\,V\\		
		\bottomrule[1.5pt]
	\end{tabular}
\end{table}

\begin{figure*}
	\centering
	\includegraphics [scale=1] {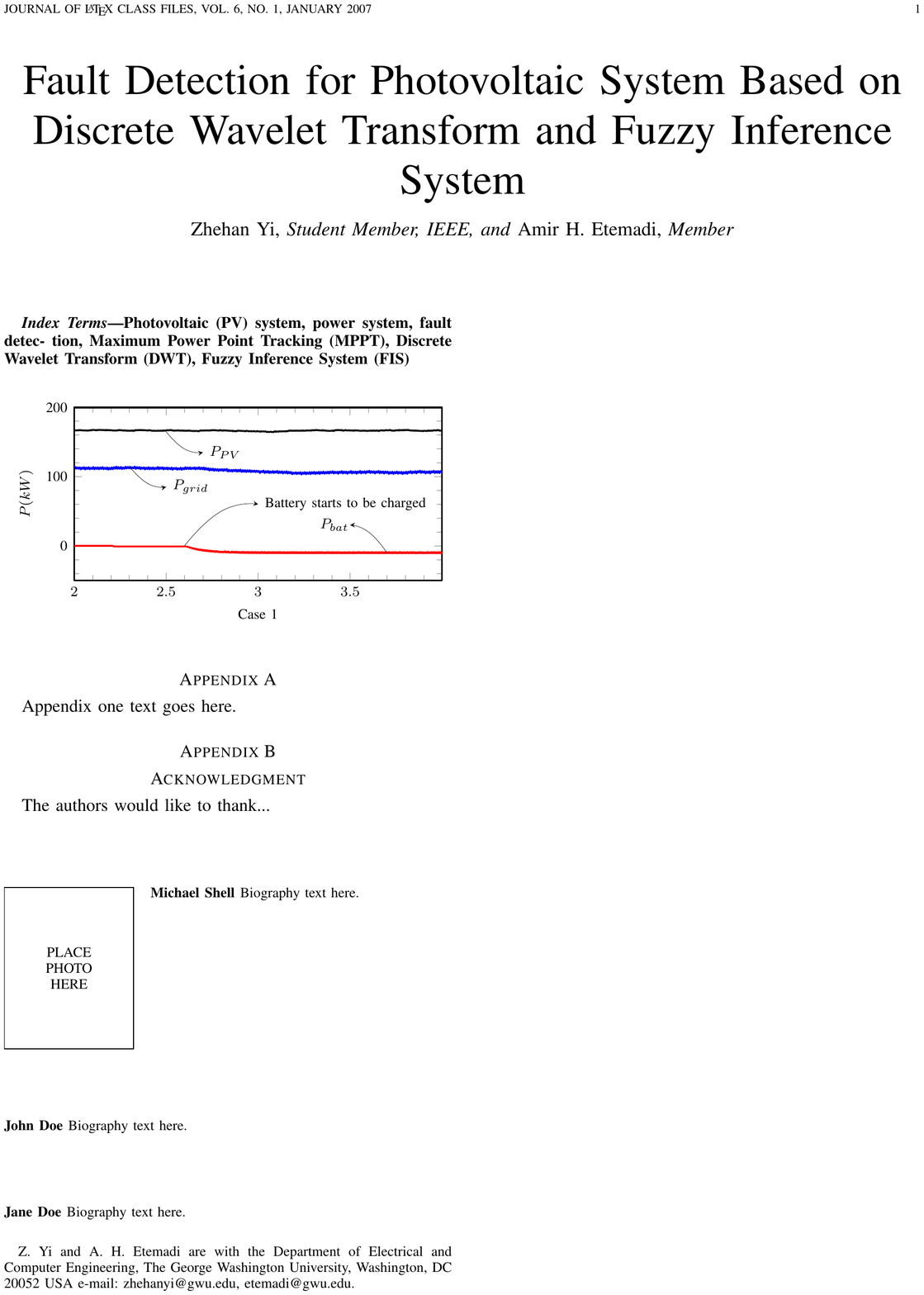}
	\includegraphics [scale=1] {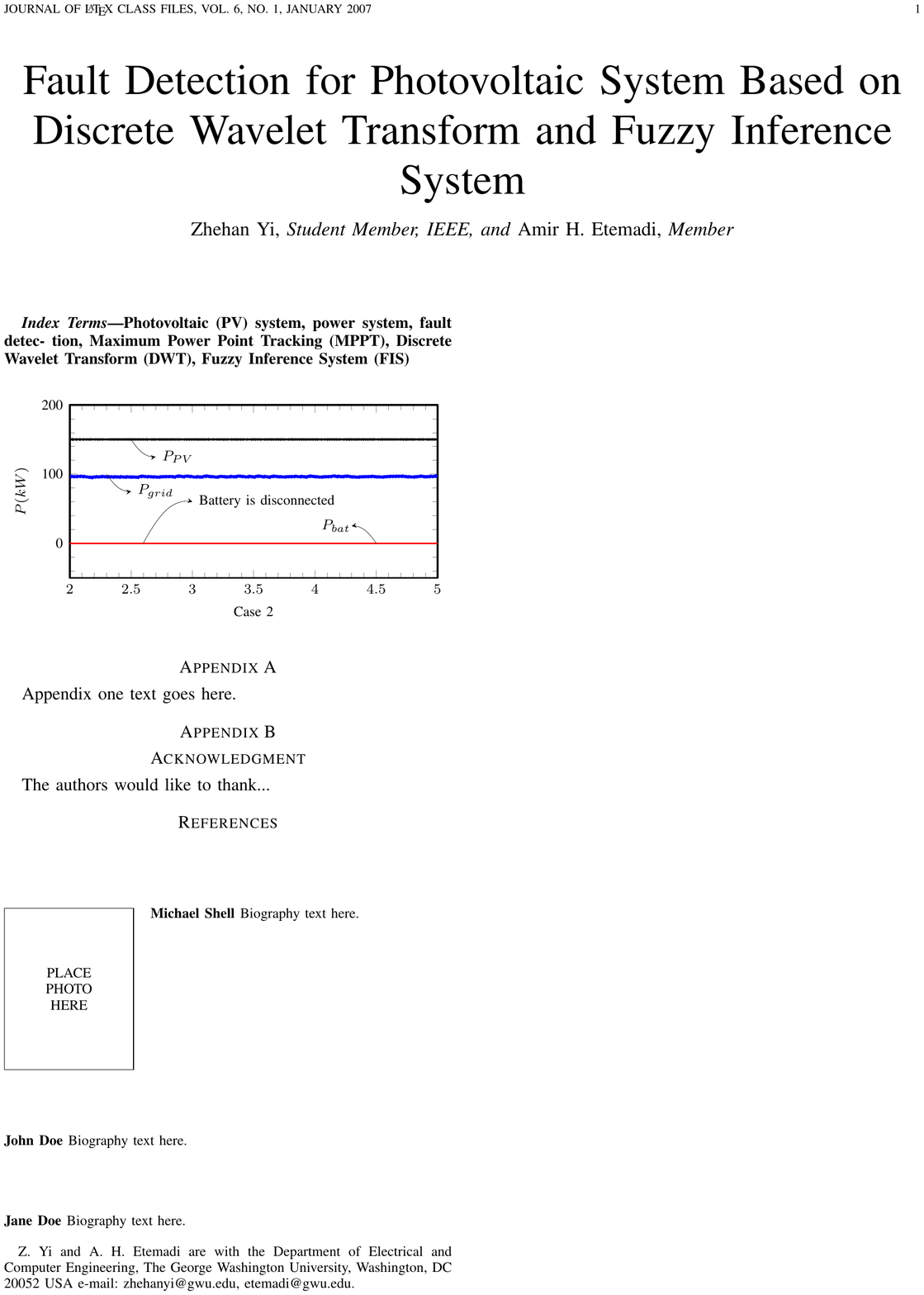}
	\includegraphics [scale=1] {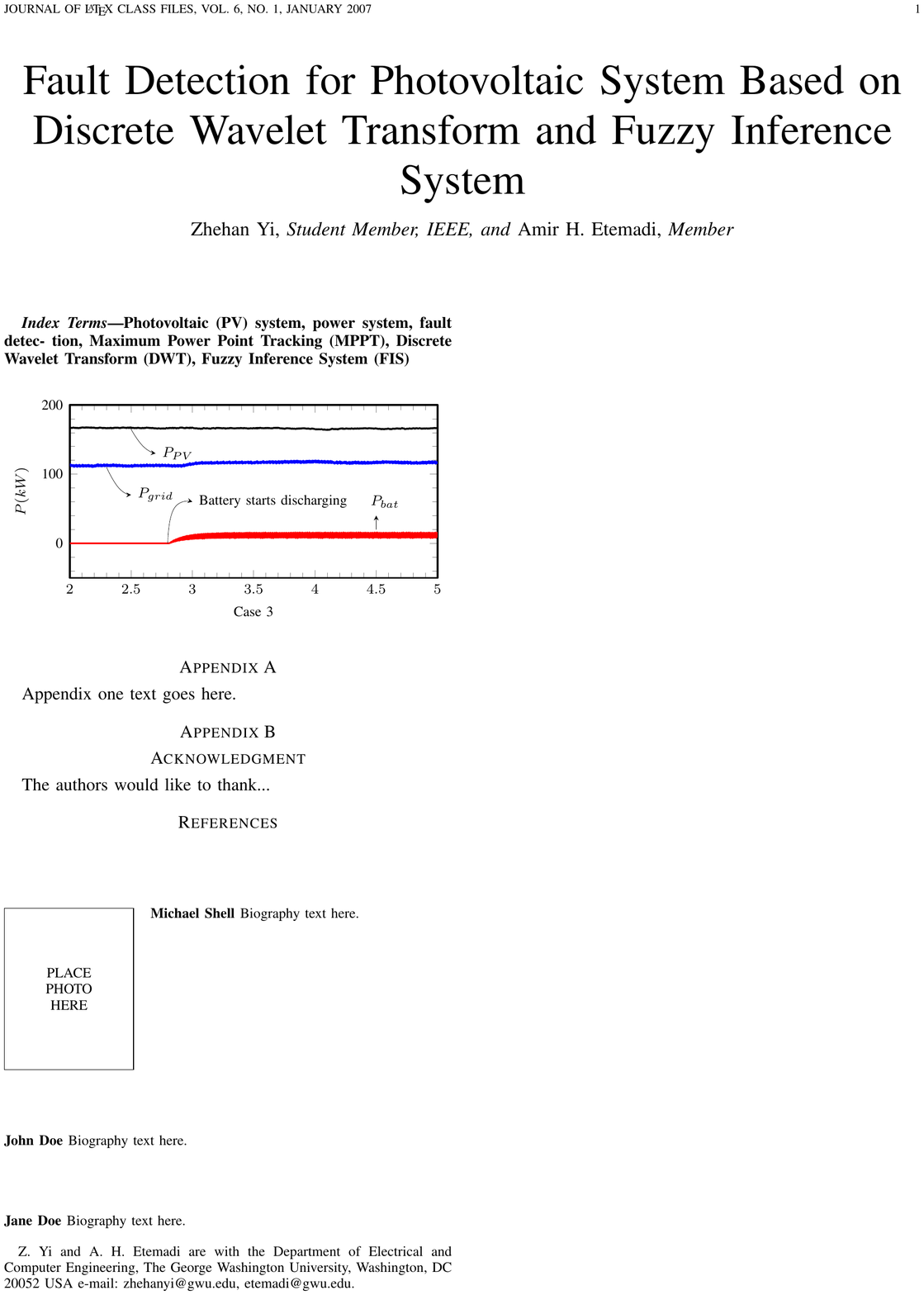}
	\includegraphics [scale=1] {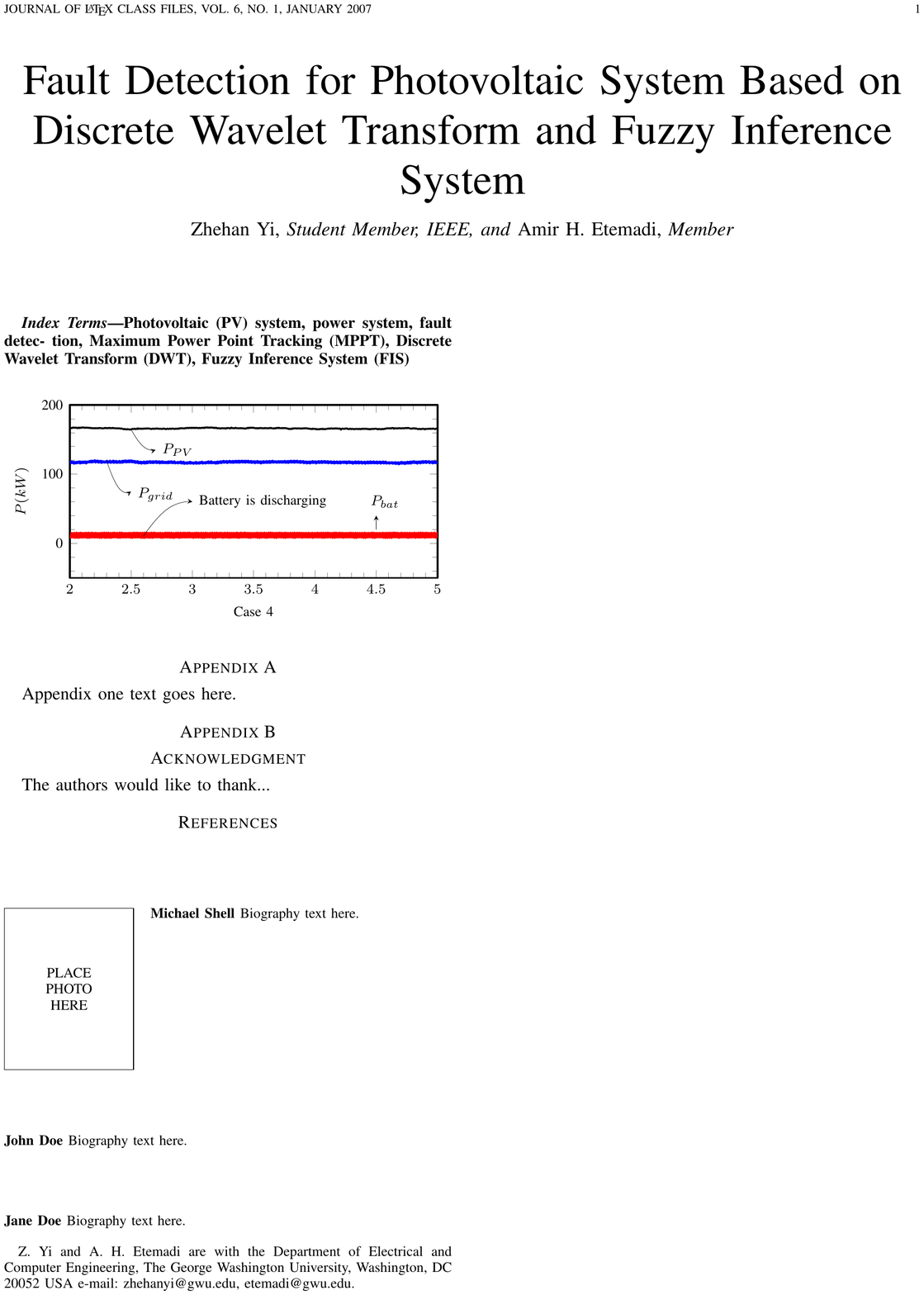}
	\includegraphics [scale=1] {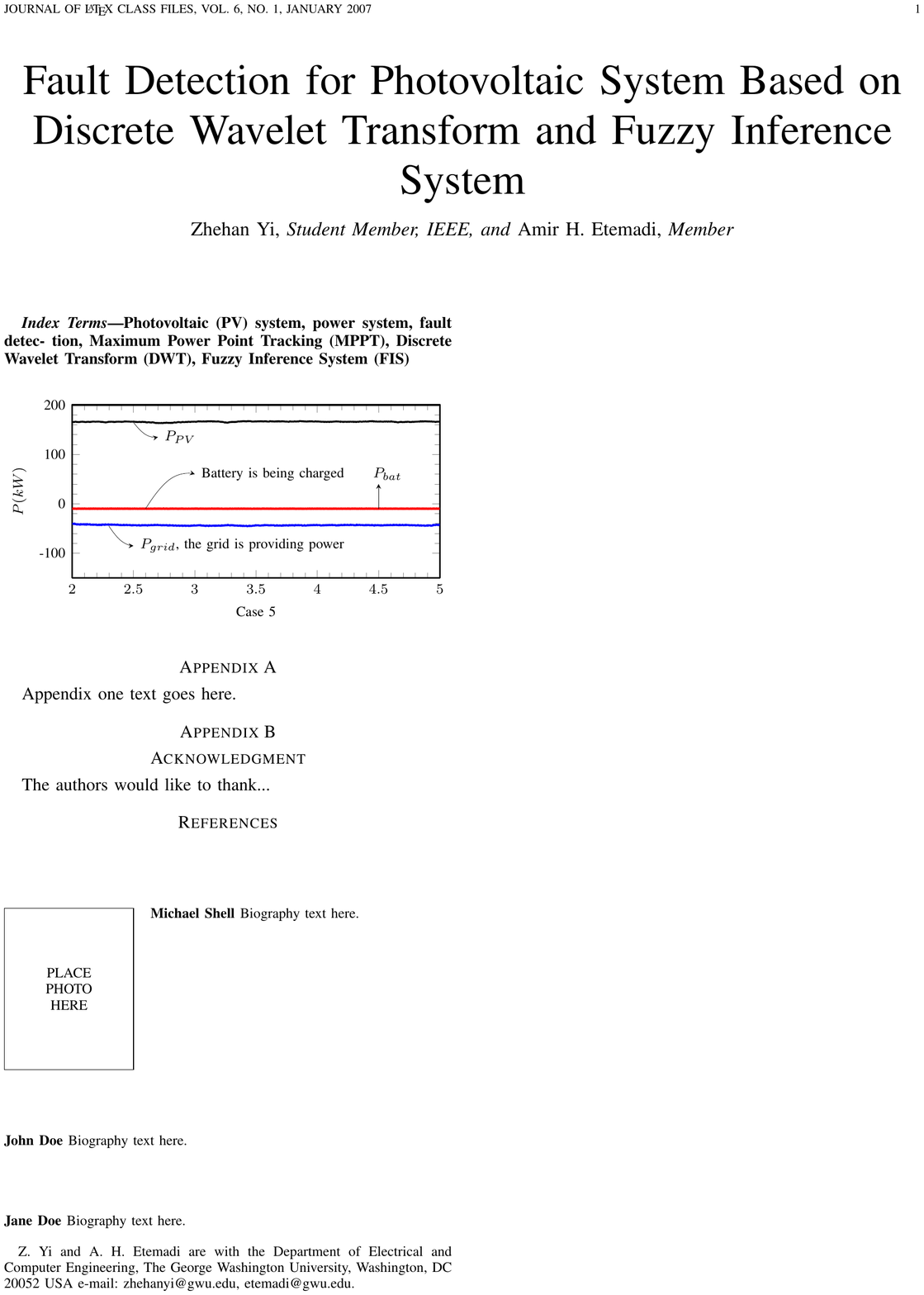}
	\caption{Power flows in the case studies (the $x$ axis represents time in second).}
	\label{case}
\end{figure*}

\subsubsection*{Case 1}

When there is excess power from the PV array after fulfilling the DC load and AC demand, while the battery is not fully charged ($P_{PV}>P_{DC-load}+P_{grid}$, $\textup{SOC} < 95\%$), the excess PV power will be transferred and stored in the battery. The power flows in the system are shown in Fig. \ref{case} - Case 1, which illustrates the balanced power in the system: the PV array is working in MPPT mode which provides $P_{PV}$ (165\,kW in Fig. \ref{case} - Case 1), while the DC load is constant at 50\,kW, $P_{grid}$ is around 105\,kW, and the battery is charged by approximately 10\,kW.

\subsubsection*{Case 2} 

When the maximum power from the PV array is greater than the sum of DC load and AC demand, and the battery is fully charged ($P_{PV}>P_{DC-load}+P_{grid}$, $\textup{SOC} \geq 95\%$), the PV array will switch from MPPT control mode to reference power control mode, where the reference is $P_{PV-ref} = P_{DC-load}+P_{grid}$, and the battery is disconnected. As is shown in Fig. \ref{case} - Case 2, $P_{PV}$ is controlled at 150\,kW, the battery is neither charging nor discharging ($P_{bat}$), the DC load consumes 50\,kW,  and the remaining PV power is feeding the grid (100\,kW).

\subsubsection*{Case 3} 

For the case where the demand from DC load and the grid is greater than the maximum PV power, and the battery is not over-discharged ($P_{PV}<P_{DC-load}+P_{grid}<P_{PV}+P_{bat}$, $\textup{SOC} \geq 20\%$), the PV array will work under MPPT control, providing the 165\,kW, and the battery will discharge by 10\,kW to compensate the difference between supply and demand (Fig. \ref{case} - Case3).

\subsubsection*{Case 4} 

If the total demand from DC load and the grid is greater than the maximum PV power and the battery power, ($P_{PV}+P_{bat}<P_{DC-load}+P_{grid}$, $\textup{SOC} \geq 20\%$), similar to Case 3, the PV array will work on MPPT mode and the battery will discharge, providing as much power to the grid as the system can (Fig. \ref{case} - Case4).

\subsubsection*{Case 5} 

In Fig. \ref{case} - Case 5, where the PV array cannot fulfill the DC load (increased from 50\,kW to 190\,kW), the grid does not request power from the DC microgrid, and the battery has no excess power ($P_{PV}<P_{DC-load}$, $\textup{SOC} \leq 20\%$), while maintain MPPT generation in the PV array (165\,kW), additional electric power can be purchased from the grid (35\,kW) to meet the demand on the DC side, and, when necessary, to charge the battery (10\,kW) for later needs.

The case studies successfully presents the satisfactory performance of the proposed power control and management system. Power flows in all circumstances mentioned above are properly balanced. In all the cases, whether the PV array is controlled under power reference or MPPT mode, demands from both the DC-load and the utility grid are reliably supplied by controlling the DC/DC converters and the DC/AC inverter, and equation (\ref{eq1}) is always maintained. Although waveforms are not shown, the DC bus voltage $V_{DC}$ is stabilized around 450\,V to ensure a stable DC power supply, and the reactive power $Q$ through inverter is controlled at 0\,Var in all these cases, regardless of any change in active power flows. Nevertheless, $V_{DC}$ can be easily changed to any reasonable value, by simply modifying the reference $V_{DC-ref}$ in Fig. \ref{inverter}. When requested, the system can provide reactive power $Q$ to the grid quickly, by setting the reference $Q_{ref}$ (Fig. \ref{inverter}).      

\section{Conclusion}
Power management in grid-connected PV-battery systems is critical to maintain a reliable power supply to load and the utility grid demands. This paper proposes a power control and management system, which is able to effectively manage the power flows in grid PV-battery microgrid systems. Power demands and supplies are successfully balanced by the control of power converts, and the reactive power is also under full monitoring and control. The proposed system regulates the DC bus voltage by controlling the inverter, such that power can be provided to feed the load reliably in spite of other changes. The DC bus voltage value is under full control. Additionally, it is more convenient and economical for DC load access, as DC/DC converters can be omitted if the rating voltage of the load meets the DC bus voltage. Case studies are carried out, and the performance of the proposed system is successfully verified.






%

\end{document}